\documentclass[aos,dvips]{arximspdf}
\usepackage{mathrsfs}
\usepackage{graphicx}

\doi{10.1214/10-AOS858}
\volume{39}
\issue{2}
\pubyear{2011}
\firstpage{1042}
\lastpage{1068}

\makeatletter

\newtheorem{theorem}{Theorem}

\newtheorem{pro}{Proposition}
\newproclaim{rem}{Remarks}

\newproclaim{example}{Example}

\newcommand{\R}{\mathbb{R}}
\newcommand{\E}{\mathbb{E}}
\newcommand{\eps}{\varepsilon}

\newcommand{\bS}{\mathbb{S}}

\newcommand{\bN}{\mathbb{N}}

\newcommand{\EE}{\mathbb{E}}
\newcommand{\bP}{\mathbb{P}}
\newcommand{\bE}{\mathbb{E}}
\newcommand{\bL}{\mathbb{L}}

\newcommand{\tep}{{t_N}}

\newcommand{\bjk}{\beta_{j\eta}}
\newcommand{\psijk}{\psi_{j\eta}}

\makeatother

\begin{document}
\begin{frontmatter}

\title{Localized spherical deconvolution}
\runtitle{Localized spherical deconvolution}

\begin{aug}
\author[A]{\fnms{G\'{e}rard} \snm{Kerkyacharian}\ead[label=e1]{kerk@math.jussieu.fr}},
\author[B]{\fnms{Thanh~Mai}~\snm{Pham~Ngoc}\corref{}\ead[label=e2]{thanh.pham\_ngoc@math.u-psud.fr}}
\and
\author[C]{\fnms{Dominique} \snm{Picard}\ead[label=e3]{picard@math.jussieu.fr}}
\runauthor{G. Kerkyacharian, T. M. Pham Ngoc and D. Picard}
\affiliation{CNRS and LPMA, Universit\'{e} Paris Sud and Universit\'{e} Paris VII}
\address[A]{G. Kerkyacharian\\
CNRS, LPMA\\
175 rue du Chevaleret\\
75013 Paris\\
France\\
\printead{e1}} 
\address[B]{T. M. Pham Ngoc\\
Laboratoire de Math\'{e}matiques,\\
\quad UMR CNRS 8628\\
Universit\'{e} Paris Sud\\
91405 Orsay Cedex\\
France\\
\printead{e2}}
\address[C]{D. Picard\\
Universit\'{e} Paris VII\\
175 rue du Chevaleret\\
75013 Paris\\
France\\
\printead{e3}}
\end{aug}

\received{\smonth{1} \syear{2010}}
\revised{\smonth{9} \syear{2010}}

%
\begin{abstract}
We provide a new algorithm for the treatment of the deconvolution
problem on the sphere which combines the traditional SVD inversion with
an appropriate thresholding technique in a well chosen new basis. We
establish upper bounds for the behavior of our procedure for any
$\bL_p$ loss. It is important to emphasize the adaptation properties
of our procedures with respect to the regularity (sparsity) of the
object to recover as well as to inhomogeneous smoothness. We also
perform a numerical study which proves that the procedure shows very
promising properties in practice as well.
\end{abstract}

%

\end{frontmatter}
%


\section{Introduction}
The spherical deconvolution problem was first proposed by Rooij and
Ruymgaart \cite{vanRooijRuymgaart} and subsequently solved in
Healy, Hendriks and Kim~\cite{HealyHendriksKim}. Kim and Koo
\cite{Kimoptimalspherical} established minimaxity for the
$\bL_2$-rate of convergence. The optimal procedures obtained there are
using orthogonal series methods associated with spherical harmonics.
One important problem arising with these procedures is their poor
local performances due to the fact that spherical harmonics are spread
all over the sphere. This explains for instance the fact that although
they are optimal in the $\bL_2$ sense, they cease to be optimal for
other losses, such as~$\bL_p$ losses for instance.

In our approach, we focus on two important points. We aim at a
procedure of estimation which is efficient from a $\bL_2$ point of
view, as well as it performs satisfactorily from a local point of view
(for other $\bL_p$ losses for instance).

Deconvolution is an inverse problem and in such there is a notable
conflict between the inversion part which in presence of noise creates
an instability reasonably handled by a Singular Value Decomposition
(SVD) approach and
the fact that the SVD basis is very rarely localized and capable of
representing local features of images,
which are especially important to recover. Our strategy is to follow
the approach started in Kerkyacharian et al.
\cite{KerkPetruPicaWiller} for the Wicksell case, Kerkyacharian et al.
\cite{5authors} for the Radon transform, which utilizes the
construction borrowed from Narcowich, Petrushev and Wald
\cite{NarcoPetruWard,NPW} of a tight frame (i.e., a redundant
family) staying sufficiently close to the SVD decomposition but which
enjoys at the same time enough localization properties to be
successfully used for statistical estimation (see, for instance, Baldi
et al. \cite{BaldiKerkMariPica,subsampling}, Pietrobon et al. \cite{notice}  for other types of applications). The construction
\cite{NPW} produces a family of functions which very much resemble
wavelets, the needlets.

To achieve the goals presented above, and especially adaptation to
different regularities and local inhomogeneous smoothness, we
essentially use a~projection method on the needlets (which enables a
stable inversion of the deconvolution, due to the closeness to the SVD
basis) with a subsequent fine tuning thresholding process.

This provides a reasonably simple algorithm with very good
performances, both from a theoretical point of view and a numerical
point of view. In effect, this new algorithm provides a much better
spatial adaptation, as well as adaptation to wider classes of
regularity. We give here upper bounds obtained by the procedure over a
large class of Besov spaces and any $\bL_p$ losses as well as
$\bL_\infty$. We find back these results in the simulation study where
the effect of the localization are highlighted for instance by a
comparison of the performances, on a bell-density example between the
procedure provided here and the SVD methods (detailed below) proving
that our quality of reconstruction of the peak is notably better.

It is important to notice that especially because we consider different
$\bL_p$ losses, we provide rates of convergence of new types attained
by our procedure, which, of course, coincide with the usual ones for
$\bL_2$ losses.

Again, the problem of choosing appropriated spaces of regularity on
the sphere in a serious question, and we decided to consider the spaces
which may be the closest to our natural intuition: those which
generalize to the sphere case the classical approximation properties of
usual regularity spaces such as H\"{o}lder spaces and include at the
same time the Sobolev regularity spaces used in Kim and Koo~\cite{Kimoptimalspherical}.

Sphere deconvolution has a vast domain of application such as medical
imaging (see Tournier et al.
\cite{TournierCalamanteGadianConnelly}) and astrophysics. Indeed, our
results are especially motivated by many recent developments in the
area of observational astrophysics.

It is a common problem in astrophysics to analyze data sets consisting
of a number of objects (such as galaxies of a particular type) or of
events (such as cosmic rays or gamma ray bursts) distributed on the
celestial sphere. In many cases, such objects trace an underlying
probability distribution $f$ on the sphere, which itself depends on the
physics which governs the production of the objects and events.

The case for instance of ultra high energy cosmic rays (UHECR)
illustrates well the type of applications covered by our results. Ultra
high energy cosmic rays are particles of unknown nature which arrive at
the earth from apparently random directions of the sky. They could
originate from long-lived relic particles from the Big Bang, about 13
billion years old. Alternatively, they could be generated by the
acceleration of standard particles, such as protons, in extremely
violent astrophysical phenomena, such as cluster shocks. They could
also originate from Active Galactic Nuclei (AGN), or from neutron stars
surrounded by extremely high magnetic fields.

Hence, in some hypotheses, the underlying probability distribution for
the directions of incidences of observed UHECRs would be a finite sum
of point-like sources---or near point like, taking into account the
deflection of the cosmic rays by magnetic fields. In other hypotheses,
the distribution could be uniform, or smooth and correlated with the
local distribution of matter in the universe. The distribution could
also be a superposition of the above. Identifying between these
hypotheses is of primordial importance for understanding the origin and
mechanism of production of UHECRs.

Of course, the observations of these events ($X_i$'s in the sequel) are
always most often perturbated by a secondary noise ($\eps_i$) which
leads to the deconvolution problem described below. Following Healy,
Hendriks and Kim \cite{HealyHendriksKim}, Kim and Koo
\cite{Kimoptimalspherical}, the spherical deconvolution problem can
be described as follows. Consider the situation where we observe
$Z_1,\ldots, Z_N,$ $N$ i.i.d. observations with
\begin{equation}\label{1}
Z_i=\eps_i X_i,
\end{equation}
where the $\eps_i$'s are i.i.d. random elements in $\mathit{SO}(3)$ (the group
of $3\times3$ rotation matrices), and the $Z_i$'s and $X_i$'s are
i.i.d. random elements of $\mathbb{S}^2$ (two-dimensional unit sphere
of $\R^3$) random elements, with $\eps_i$ and $X_i$ assumed to be
independent. We suppose that the distributions of $X$ and $Z$ are
absolutely continuous with respect to the uniform probability measure
on $\mathbb{S}^2$, and that the distribution of $\varepsilon$ is
absolutely continuous with respect to the Haar measure of $\mathit{SO}(3)$. We
will denote the densities of resp. $Z, X$ and $\varepsilon$ resp.
$f_Z, f_X,  f_\eps$.

Then
%
\begin{equation}
f_Z=f_\eps*f_X,
\end{equation}
where $*$ denotes convolution and is defined below. In the sequel,
$f_X$ will be denoted by $f$ to emphasize the fact that it is the
object to recover.

The following paragraph recalls the necessary definitions. It is
largely inspired by Kim and Koo \cite{Kimoptimalspherical}
and Healy,
Hendriks and Kim \cite{HealyHendriksKim}.

\section{Some preliminaries about harmonic analysis on $\mathit{SO}(3)$ and $\mathbb S^2$}

We will provide a brief overview of Fourier analysis on $\mathit{SO}(3)$ and
$\mathbb{S}^2$. Most of the material can be found in an expanded form
in Vilenkin \cite{Vilenkin}, Talman~\cite{Talman}, Terras
\cite{Terras}, Kim and Koo \cite{Kimoptimalspherical},
and Healy,
Hendriks and Kim \cite{HealyHendriksKim}. Let
\[
u(\phi)=\pmatrix{
\cos\phi& -\sin\phi& 0 \cr
\sin\phi& \cos\phi& 0 \cr
0 & 0 & 1}, \qquad a(\theta)=\pmatrix{
\cos\theta& 0 & \sin\theta\cr
0 & 1 & 0 \cr
-\sin\theta& 0 & \cos\theta
},
\]
where, $\phi\in[0,2\pi),\theta\in[0,\pi)$. It is well known
that any
rotation matrix can be decomposed as a product of three elemental
rotations, one about the $z$-axis first by an angle $\psi$, followed by
a rotation about the $y$-axis by an angle $\theta$, and finally by
another rotation again about the $z$-axis by an angle~$\phi$. Indeed,
the well-known Euler--Angle decomposition says that any $g \in \mathit{SO}(3)$
can almost surely be uniquely represented by three angles
$(\phi,\theta, \psi)$, with the following formula (see Healy,
Hendriks and Kim \cite{HealyHendriksKim} for details):
%
\begin{equation}\label{decompositionSO3}
g=u(\phi)a(\theta)u(\psi),
\end{equation}
where $\phi\in[0,2\pi),\theta\in[0,\pi), \psi\in[0,2\pi)$.
Consider the functions, known as the rotational harmonics,
%
\begin{equation}\label{harmo_rotationnelle}
D^l_{mn}(\phi,\theta, \psi)=e^{-i(m\phi+n\psi
)}P^l_{mn}(\cos\theta),
\end{equation}
where the associated Legendre functions $P^l_{mn}$ for $-l\le m,n\leq
l,l=0,1,\ldots,$ are fully described in Vilenkin
\cite{Vilenkin}. The functions $D^l_{mn}$ for $-l\le m,n\leq
l,l=0,1,\ldots,$ are the eigenfunctions of the Laplace Beltrami
operator on $\mathit{SO}(3)$, hence, $\sqrt{2l+1} D^l_{mn}, -l\le m,n\leq
l,l=0,1,\ldots$ is a complete orthonormal basis for
$\mathbb{L}_2(\mathit{SO}(3))$ with respect to the probability Haar measure. In
addition, if we define the $(2l+1)\times(2l+1)$ matrices by
%
\begin{equation}\label{3}
D^l (g)=[D^l_{mn}
(g)],
\end{equation}
where for $-l\le m,n\leq l,l=0,1,\ldots,$ and $g \in \mathit{SO}(3)$, they
constitute the collection of inequivalent irreducible representations
of $\mathit{SO}(3)$ (for further details, see Vilenkin \cite{Vilenkin}).

Hence, for $f \in\mathbb{L}_2(\mathit{SO}(3))$, we define the rotational
Fourier transform on $\mathit{SO}(3)$ by
%
\begin{equation}\label{rotationalFourierTransform}
\hat f_{mn}^l=\int_{\mathit{SO}(3)}f(g)D^l_{mn}(g)\,dg,
\end{equation}
where $dg$ is the probability Haar measure on $\mathit{SO}(3)$ and we define
the following matrix of dimension $(2l+1)\times(2l+1)$
\[
\hat f^l=[\hat f_{mn}^l]_{-l\le m,n\leq l},\qquad l=0,1,\ldots.
\]
The rotational inversion can be obtained by
\begin{eqnarray}\label{4}
f(g)&=&\sum_l\sum_{-l\le m,n\leq l}\hat f_{mn}^l{\overline
{D^l_{mn}(g)}}\nonumber\\ [-8pt]\\ [-8pt]
&=&\sum_l\sum_{-l\le m,n\leq l}\hat f_{mn}^l{
D^l_{mn}(g^{-1})},\nonumber
\end{eqnarray}
(\ref{4}) is to be understood in $\mathbb{L}_2$-sense although with
additional smoothness conditions, it can hold pointwise.

A parallel spherical Fourier analysis is available on $\mathbb{S}^2$.
Any point on $\mathbb{S}^2$ can be represented by
\[
\omega=(\cos\phi\sin\theta,\sin\phi\sin\theta,\cos\theta)^t,
\]
with, $\phi\in[0,2\pi),\theta\in[0,\pi)$. We also define the
functions
%
\begin{equation}\label{5}
Y^l_m(\omega)=Y^l_m(\theta,\phi)=(-1)^m\sqrt{\frac{(2l+1)}{4\pi
}\frac{(l-m)!}{(l+m)!}}%
P^l_m(\cos\theta)e^{im\phi}
\end{equation}
for $-l\le m\leq l,l=0,1,\ldots,$ $\phi\in[0,2\pi),\theta
\in
[0,\pi)$ and where $P^{l}_{m}(\cos
\theta)$ are the Legendre functions. 

The functions $Y^{l}_{m}$ obey
%
\begin{equation}\label{harmonique_negative}
Y^{l}_{-m}(\theta,\phi)={(-1)}^m\overline{{Y}^{l}_{m}(\theta,\phi)}.
\end{equation}
The set $\{Y^l_m, -l\le m\leq l,l=0,1,\ldots\}$ is forming an
orthonormal basis of $\mathbb{L}_2(\mathbb{S}^2)$, generally referred
to as the spherical harmonic basis.

Again, as above, for $f \in\mathbb{L}_2(\mathbb{S}^2)$, we define the
spherical Fourier transform on $\mathbb{S}^2$ by
%
\begin{equation}\label{7}
\hat f_{m}^l=\int_{\mathbb{S}^2}f(\omega){\overline{
Y^l_{m}(\omega)}}\,d\omega,
\end{equation}
where $d\omega$ is the uniform probability measure on the sphere
$\mathbb{S}^2$.
The spherical inversion can be obtained by
\begin{equation}\label{7}
f(\omega)=\sum_l\sum_{-l\le m \leq l}\hat f_{m}^l{Y^l_{m}(\omega)}.
\end{equation}
The bases detailed above are important because they realize a singular
value decomposition of the convolution operator created by our model.
In effect, we define for $f_\eps\in\mathbb{L}_2(\mathit{SO}(3)),  f\in
\mathbb{L}_2(\mathbb{S}^2)$ the convolution by the following formula:
\[
f_\eps*f(\omega)=\int_{\mathit{SO}(3)}f_\eps(u)f(u^{-1}\omega)\,du,
\]
and we have for all $-l\le m\leq l,l=0,1,\ldots,$
%
\begin{equation}\label{convprod}
(\widehat{f_\eps*f})_m^l=\sum_{n=-l}^l\hat f^l_{\eps,mn}\hat
f^l_n:=(\hat f_\eps^l\hat f^l)_m .
\end{equation}

\subsection{The SVD method}
The singular value method (see Healy, Hendriks and Kim
\cite{HealyHendriksKim} and Kim and Koo
\cite{Kimoptimalspherical}) consists in expanding $f$ in the
spherical harmonics basis $Y^l_{m}$ and estimating the spherical
Fourier coefficients using the formula above~(\ref{convprod}). We get
the following estimator of the spherical Fourier transform of~$f$:
\begin{eqnarray}\label{empcoef}
\hat f^{l,N}_{m}&:=&\frac1N\sum_{j=1}^N\sum_{n=-l}^l\hat
f^{l}_{\eps^{-1},mn}\bar{Y}^l_n(Z_j),\nonumber\\ [-8pt]\\ [-8pt]
\hat f^{l}_{\eps^{-1}}&:=&(\hat f^{l}_{\eps})^{-1}\nonumber,
\end{eqnarray}
provided, of course, that these inverse matrices exist, and then the
estimator of the distribution $f$ is
%
\begin{equation}
f^{N}(\omega)=\sum_{l=0}^{\tilde{N}}\sum_{m=-l}^l
\hat
f^{l,N}_{m}Y^l_m(\omega),
\end{equation}
where $\tilde{N}$ is depending on the number of observations and has to
be properly selected.

\section{Needlet construction}\label{sec3}

This construction is due to Narcowich et al.~\cite{NarcoPetruWard}.
Its aim is essentially to build a very well
localized tight frame constructed using spherical harmonics, as
discussed below. It was recently extended to more general Euclidean
settings with fruitful statistical applications (see Kerkyacharian et
al.~\cite{KerkPetruPicaWiller}, Baldi et al. \cite{BaldiKerkMariPica,subsampling}, Pietrobon
et al. \cite{notice}). As
described above, we have the following decomposition:
%
\begin{equation}\label{decomp}
\mathbb{L}_2(\mathbb{S}^2)=\bigoplus_{l=0}^\infty\mathbb{H}_l,
\end{equation}
where $\mathbb{H}_l$ is the space spanned by $\{Y^l_m, -l\le m\leq
l\}$ of spherical harmonics of $\mathbb{S}^2$, of degree $l$ (which
dimension is $2l+1$).

The orthogonal projector on $\mathbb{H}_l$ can be written using the
following kernel operator:
%
\begin{equation}
\forall\! f\in\mathbb{L}_2(\mathbb{S}^2)\qquad
P_{\mathbb{H}_l}f(x)=\int_{\mathbb{S}^2}L_l(\langle x,y\rangle)
f(y)\,dy,
\end{equation}
where,
\[
L_{l}( x,y )=\sum_{m=-l}^{l} Y^{l}_{m}(x) \overline
{Y^{l}_{m}(y)}=L_{l}(\langle x,y \rangle),
\]
and where $\langle x,y\rangle$ is the
standard scalar product of $\mathbb{R}^{3}$, and $L_l$ is the Legendre
polynomial of degree $l$, defined on $[-1,+1]$ and verifying
%
\begin{equation}\label{lnorm}
\int_{-1}^1L_l(t)L_k(t) =\frac{2l+1}{8\pi}\delta_{l,k},
\end{equation}
where $\delta_{l,k}$ is the Kronecker symbol.

Let us point out the following reproducing property of the projection
operator:
%
\begin{equation}\label{auto}
\int_{\bS^2} L_{l}(\langle x,y\rangle)L_{k}(\langle
y,z\rangle)\,dy=\delta_{l,k}L_{l}(\langle x,z\rangle).
\end{equation}

The following construction is based on two fundamental steps:
Littlewood--Paley decomposition and discretization, which are summarized
in the two following subsections.

\subsection{Littlewood--Paley decomposition}\label{ssec-littlewood}

Let $\phi$ be a $C^{\infty}$ function on $\mathbb R$, symmetric and
decreasing on $\mathbb R^+$ supported in $|\xi|\leq1,$ such that
$1\geq\phi(\xi)\geq0$ and $\phi(\xi)=1$ if $|\xi|\leq
\frac{1}{2} $.
\[
b^{2}(\xi)=\phi\biggl(\frac\xi2\biggr)-\phi(\xi)\geq0,
\]
so that
%
\begin{equation} \label{1}
\forall|\xi|\geq1\qquad\sum_{j\ge0}b^{2}\biggl(\frac{\xi
}{2^{j}}\biggr)=1.
\end{equation}

Remark that $b(\xi)\not=0$ only if $\tfrac12\le|\xi|\le2$. Let us now
define the operator
$\Lambda_{j}=\sum_{l\geq0}b^{2}(\tfrac{l}{2^{j}})L_{l} $ and the
associated kernel
\[
\Lambda_{j}(x,y)=\sum_{l\geq0}b^{2}\biggl(\frac{l}{2^{j}}\biggr)L_{l}(\langle
x,y\rangle)=%
\sum_{2^{j-1}<l<2^{j+1}}b^{2}\biggl(\frac{l}{2^{j}}\biggr)L_{l}(\langle
x,y\rangle).
\]
We obviously have
%
\begin{equation}\label{rep}
\forall\! f\in{\bL_2(\mathbb{S}^2)}\qquad f=\lim
_{J\rightarrow
\infty}L_{0}(f)+\sum_{j=0}^{J}\Lambda_{j}(f)
\end{equation}
and if $M_{j}(x,y)=\sum_{l\geq0}b(\frac{l}{2^{j}})L_{l}(\langle
x,y\rangle)$, then
%
\begin{equation}\label{sqrt}
\Lambda_{j}(x,y)=\int M_{j}(x,z)M_{j}(z,y)\,dz.
\end{equation}

\subsection{Discretization and localization properties}

Let us define
\[
\mathscr{P}_{l}=\bigoplus_{m=0}^{l}\mathbb{H}_{m},
\]
the space spanned by the spherical harmonics of of degree less than
$l$.

The following quadrature formula is true: for all $l\in\mathbb{N}$
there exists a finite subset $\mathscr{X}_{l}$ of $S^2$ and positive
real numbers $ \lambda_{\eta}>0$, indexed by the elements $\eta$ of
${}\mathscr{X}_{l},$ such that
%
\begin{equation}\label{quadr}
\forall\! f\in\mathscr{P}_{l}\qquad\int_{\mathbb{S}^2} f(x)\,
dx=\sum_{\eta\in\mathscr{X}%
_{l}}\lambda_{\eta}f(\eta).
\end{equation}

Then the operator $M_{j}$ defined in the subsection above is such that
\[
z\mapsto M_{j}(x,z)\in\mathscr{P}_{[2^{j+1}]},
\]
so that
\[
z\mapsto M_{j}(x,z)M_{j}(z,y)\in\mathscr{P}_{[2^{j+2}]}
\]
and we can write
\[
\Lambda_{j}(x,y)=\int M_{j}(x,z)M_{j}(z,y)\,dz=\sum_{\eta\in\mathscr
{X}_{[2^{j+2}]}}\lambda_{\eta}M_{j}(x,\eta)M_{j}(\eta,y).
\]
This implies
\begin{eqnarray*}
\Lambda_{j}f(x) & =&\int\Lambda_{j}(x,y)f(y)\,dy=\int\sum_{\eta\in%
\mathscr{X}_{[2^{j+2}]}}\lambda_{\eta}M_{j}(x,\eta)M_{j}(\eta,y)f(y)\,dy\\
& =&\sum_{\eta\in\mathscr{X}_{[2^{j+2}]}}\sqrt{\lambda_{\eta}}%
M_{j}(x,\eta)\int{\sqrt{\lambda_{\eta}}M_{j}(y,\eta)}f(y)\,dy.
\end{eqnarray*}
We denote
\[
\mathscr{X}_{[2^{j+2}]}=\mathscr{Z}_j,\qquad\psi_{j,\eta}(x):=
\sqrt{\lambda_{\eta}}M_{j}(x,\eta)\qquad\mbox{for }\eta\in\mathscr{Z}_j.
\]

It can also be proved that the set of cubature points $\mathscr{X}_{l}$
can be chosen so that
%
\begin{equation}\label{card}
\frac1c  2^{2j}\le\# \mathscr{Z}_j\le c 2^{2j},\qquad \frac
1c  2^{2j}\le\lambda_{\eta}\le c 2^{2j}
\end{equation}
for some $c>0$. It holds, using (\ref{rep}):
\[
f = L_{0}(f) + \sum_{j}\sum_{\eta\in\mathscr{Z}_j} \langle
f,\psi_{j,\eta}\rangle_{\bL_2(\bS^2)}\psi_{j,\eta}.
\]

The main result of Narcowich et al. \cite{NarcoPetruWard} is the
following localization property of the
$\psi_{j,\eta}$, called needlets: for any $k \in\bN$ there exists
a constant $c_{k}$ such that, for every $\xi\in
{\mathbb S}^2$,
%
\begin{equation}\label{localized}
|\psi_{j,\eta}(\xi)|\leq\frac{c_{k}2^{j}}{(1+2^{j}d(\eta
,\xi))^{k}},
\end{equation}
where $d$ is the natural geodesic distance on the sphere ($d(\xi,\eta)=
\arccos\langle\eta,\xi\rangle$). In other words, needlets are almost
exponentially localized around their associated cubature point, which
motivates their name.

A major consequence of this localization property can be summarized in
the following properties which will play an essential role in the
sequel. Their proof can be found in \cite{NarcoPetruWard,NPW},
also in \cite{kyoto}.

For any $1\le p<\infty$, there exist positive constants $c_p$, $C_p$,
$c$, $C$ and $D_p$ such that
\begin{eqnarray}\label{lp}
c_p2^{2j(p/2-1)}&\le&\|\psijk\|_p^p \le C_p2^{2j(p/2-1)},\\\label{linfty}
c2^{j} &\le& \|\psijk\|_\infty\le C2^{j},\\\label{inequality}
\biggl\|\sum_{\eta\in\mathscr{Z}_j}\lambda_\eta
\psi_{j,\eta}\biggr\|_\pi &\le& c 2^{2j(1/2-1/\pi)}\Biggl(\sum_{\eta\in\mathscr{Z}_j}
|\lambda_\eta|^\pi\Biggr)^{1/\pi},\\\label{item3}
\biggl\|\sum_{\eta\in\mathscr Z_j}\lambda_\eta\psijk\biggr\|
_p^p &\le& D_p\sum_{\eta\in\mathscr Z_j}|\lambda_\eta|^p\|\psijk\|_p^p, \\\label{meyer}
\biggl\|\sum_{\eta\in\mathscr Z_j} u_{\eta}\psijk\biggr\|_\infty
&\le& C\sup_{\eta\in\mathscr Z_j}|u_{\eta}|2^{{j}}.
\end{eqnarray}

To conclude this section, let us give a graphic representation of a
spherical needlet in the spherical coordinates in order to illustrate
the above theory. In the following graphic, we chose $j=3$ and
$\eta=250$.

\subsection{Besov spaces on the sphere}

The problem of choosing appropriated spaces of regularity on the sphere
in a serious question, and we decided to consider the spaces which may
be the closest to our natural intuition: those which generalize to the
sphere case the classical approximation properties used to define for
instance Sobolev spaces. In this section, we summarize the main
properties of Besov spaces which will be used in the sequel, as
established in \cite{NarcoPetruWard}.

\begin{figure}

\includegraphics{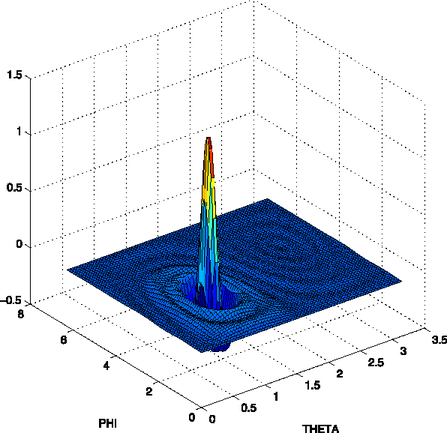}

\caption{A spherical needlet.}\label{fig1}
\end{figure}

Let $f\dvtx\mathbb{S}^2\to\R$ be a measurable function. We define
\[
E_k(f,\pi)=\inf_{P\in\mathscr{P}_k}\Vert f-P\Vert_\pi,
\]
the $\mathbb{L}_\pi$ distance between $f$ and the space
$\mathscr{P}_k$ of spherical harmonics of degree less than $k$.
The Besov space $B^s_{\pi, r}$ is defined as the space of functions
such that
\[
f\in\mathbb{L}_\pi\quad\mbox{and}\quad \Biggl(\sum_{k=0}^\infty(k^s
E_k(f,\pi))^r\frac1k\Biggr)^{1/r}<+\infty.
\]
Remarking that $k\to E_k(f,\pi)$ is decreasing, by a standard
condensation argument this is equivalent to
\[
f\in\mathbb{L}_\pi\quad\mbox{and}\quad \Biggl(\sum_{j=0}^\infty(2^{js}
E_{2^j}(f,\pi))^r\Biggr)^{1/r}<+\infty,
\]
and the following theorem states that as it is the case for Besov
spaces in $\mathbb{R}^d$, the needlet coefficients are good indicators
of the regularity and in fact Besov spaces of $\mathbb{S}^2$ are Besov
bodies, when expressed using needlet expansion. We provide in Figure \ref{fig1} a graphical
representation of a needlet which shows its localization around the
associated cubature point.
\begin{theorem} Let $1\le\pi\le+\infty$, $s>0$, $0\le r\le+\infty$. Let $f$ be a
measurable function and define
\[
\langle f,\psi_{j,\eta}\rangle=\int_{\mathbb{S}^2}f(x)
\psi_{j,\eta}(x)\,dx\stackrel{\mathit{def}.}{=}\beta
_{j,\eta},
\]
provided the integrals exist. Then $f$ belongs to $ B^s_{\pi, r}$ if
and only if, for every $j=1,2,\ldots,$
\[
\biggl(\sum_{\eta\in\mathscr{X}_j}
(|\beta_{j,\eta}|\Vert\psi_{j,\eta}\Vert_\pi)^\pi\biggr)^{1/\pi}
=2^{-js}\delta_j,
\]
where $(\delta_j)_j\in\ell_r$.
\end{theorem}

As has been shown above,
\[
c2^{2j(1/2-1/\pi)}\le\Vert\psi_{j,\eta}\Vert_\pi\le C
2^{2j(1/2-1/\pi)}
\]
for some positive constants $c,C$, the Besov space $B^s_{\pi, r}$ turns
out to be a~Banach space associated to the norm
%
\begin{equation}\label{bnorm}
\|f\|_{B^s_{\pi, r}} := \bigl\|\bigl(2^{j[s+2(1/2-
1/\pi)]}\|(\beta_{j\eta})_{\eta\in\mathscr Z_j} \|_{\ell_\pi
}\bigr)_{j\ge
0}\bigr\|_{\ell_r}<\infty,
\end{equation}
and using standard arguments (reducing to comparisons of $l_q$ norms), it
is easy to prove the following embeddings:
%
\begin{eqnarray}\label{embedding}
B^s_{\pi,r}&\subset&B^s_{p,r}\qquad\mbox{for } p\leq\pi,\nonumber\\ [-8pt]\\ [-8pt]
B^s_{\pi, r}&\subset& B^{s-2(1/\pi-1/p)}_{p,r}\qquad\mbox{for }\pi\leq p
\mbox{ and } s>2\biggl(\frac1\pi-\frac1p\biggr).\nonumber
\end{eqnarray}
Moreover, it is also true that for $s>\frac2\pi$, if $f$ belongs to
$B^s_{\pi,r}$, then it is continuous, and as a consequence bounded.

In the sequel, we shall denote by $B^s_{\pi,r}(M)$ the ball of radius
$M$ of the Besov space $B^s_{\pi,r}$.

\section{Needlet algorithm: Thresholding the needlet
coefficients}\label{thresholding}

The first step is to construct a needlet system (frame) $\{\psijk\dvtx
\eta\in\mathscr{Z}_j, j\ge-1\}$ as described in Section \ref{sec3}.

The needlet decomposition of any $f\in\bL_2(\mathbb{S}^2)$ takes the
form
\[
f=\sum_j\sum_{\eta\in\mathscr{Z}_j} (f,\psi_{j\eta})_{\bL
_2(\mathbb{S}^2)}\psi_{j\eta}.
\]
Using Parseval's identity, we have
$\beta_{j\eta}=(f,\psi_{j\eta})_{\bL_2(\mathbb{S}^2)}=
\sum_{lm} \hat{f}^l_m\psi_{j\eta}^{lm}$
with $\hat{f}^l_m=(f,Y^l_m)$ and
$\psi_{j\eta}^{lm}=(\psi_{j\eta},Y^l_m)$.

Thus,
%
\begin{equation}\label{coeff_needlet_estim}
\hat\beta_{j\eta}=\sum_{lm}\hat f^{l,N}_{m}\psi_{j\eta}^{lm},
\end{equation}
is an unbiased estimate of $\beta_{j\eta}$. We recall that $\hat
f^{l,N}_{m}$ has been defined in (\ref{empcoef}).
It is worthwhile pointing out that the SVD-estimate of the Fourier
coefficient $\hat f^{l,N}_{m}$ appears in the expression of the
estimate $\hat\beta_{j\eta}$. This underlines that the Needlet
dictionary does not depart too much from the Fourier basis and hence
benefits from the inversion property while being very well localized.

Notice that from the needlet construction (see the previous section) it
follows that the sum above is finite. More precisely,
$\psi_{j\eta}^{lm} \not= 0$ only for $2^{j-1}<l<2^{j+1}$.

Let us consider the following estimate of $f$:
%
\begin{equation}\label{estimateur_de_f}
\hat f= \sum_{j=-1}^{{J}}\sum
_{\eta\in\mathscr{Z}_j} t(\hat\beta_{j\eta})\psijk,
\end{equation}
where $t$ is a thresholding operator defined by
\begin{eqnarray}\label{def-t-N}
t(\hat\beta_{j\eta})
&=& \hat\beta_{j\eta}I\{|\hat\beta_{j\eta}|\ge{\kappa}
t_N{|\sigma_j|}\}\qquad\mbox{with} \\\label{tn}
t_N&=&\sqrt{\frac{\log N} N} ,\\\label{sigmaj}
\sigma_{j}^{2}&=&A\sum_{ln}\biggl|\sum_{m}\psi_{j\eta}^{lm}\hat
{f}^l_{\eps^{-1}mn}\biggr|^2.
\end{eqnarray}

Here $\kappa$ is a tuning parameter of the method which will be
properly selected later on. $A$ is chosen such that $\|f_Z\|_\infty\le
A$. The choice of these parameters will be discussed later. Notice that
the thresholding depends on the resolution level $j$ through the
constant $\sigma_j$. As usual in inverse problems, the upper level of
details $J$ will be chosen depending of the degree of ill-posedness. It
is precisely defined in Theorem \ref{upperbound}.

\subsection{Performances of the procedure}
The following theorem considers the case of a $\bL_p$ loss with $1\le
p<\infty$. The case $p=\infty$ is studied in Theorem \ref{upperinfty}.

\begin{theorem}\label{upperbound}
Let $1\le p<\infty$, $\nu>0$, and let us assume that
%
\begin{equation}\label{blue}
\sigma_{j}^{2} :=A\sum_{ln}\biggl|\sum_m\psi
_{j\eta}^{lm}\hat{f}^l_{\eps^{-1}mn}\biggr|^2\le C 2^{2j\nu}\qquad\forall
 j \ge0.
\end{equation}
$A$ is chosen such that $\|f_Z\|_\infty\le A$. Let us choose $\kappa$
such that $\kappa\ge\sqrt{3\pi}A$ and $ \sqrt{3\pi}A\kappa>\max\{
8p,
2p+1\}$. Let us also take $2^J= \tau[t_N]^{-1/(\nu+1)}$ with
$t_N$ as in $(\ref{tn})$ and $\tau$ is a positive constant.
Then if $\pi\ge1, $ $s > 2/\pi$, $r\ge1$ $($with the restriction
$r\le\pi$ if $s = (\nu+1)(\frac{p}{\pi} - 1))$, there exists a
constant $T$ such that
%
\begin{equation}\label{rate}
\sup_{f \in B^s_{\pi,r}(M)}\E\|\hat{ f} -f \|_p^p
\leq T(\log(N))^{p-1} \bigl[N^{-1/2} \sqrt{\log(N)}\bigr]^{\mu p},
\end{equation}
where
\begin{eqnarray*}
\mu&=&\frac{s}{s+\nu+1},\qquad  \mbox{if } s \geq
(\nu+1)\biggl(\frac{p}{\pi} - 1\biggr),\\
mu&=&\frac{s-2/\pi+2/p}{s+\nu-2/\pi+1},\qquad
\mbox{if } \frac{2}{\pi} < s < (\nu+1)\biggl(\frac{p}{\pi} - 1\biggr).
\end{eqnarray*}
\end{theorem}

The following theorem considers the case of $L_\infty$ norm loss.

\begin{theorem}\label{upperinfty}
For $\nu>0$, let us assume that there exist two constants $C$, $C'$
such that
%
\begin{equation}\label{blue'}
C' 2^{2j\nu}\le\sigma_{j}^{2} :=A\sum
_{ln}\biggl|\sum_m\psi_{j\eta}^{lm}\hat{f}^l_{\eps^{-1}mn}\biggr|^2\le C
2^{2j\nu}\qquad \forall j \ge0.
\end{equation}
$A$ is chosen such that $\|f_Z\|_\infty\le A$. Let us choose $\kappa$
such that $\kappa\ge\sqrt{3\pi}A$ and $ \sqrt{3\pi}A\kappa>16$.
Let us
also take $2^J= \tau[t_N]^{-1/(\nu+1)}$ with $t_N$ as in
$(\ref{tn})$ and $\tau$ is a positive constant.
Then
if $\pi\ge1, $ $s > 2/\pi$, $r\ge1$ , there exists a constant $T$
such that
%
\begin{equation}\label{rate}
\sup_{f \in B^s_{\pi,r}(M)}\E\|\hat{ f} -f \|_\infty
\leq T \log(N)[t_N]^{\mu' },
\end{equation}
where
\[
\mu' =\frac{s-2/\pi}{s+\nu-2/\pi+1}.
\]
\end{theorem}

The proof of these theorems is given in Section \ref{proof}.

\begin{rem*}

\begin{enumerate}
\item In Theorem \ref{upperbound}, the rates of convergence obtained
for larger $s$ are usually referred to as the dense case, whereas the
other case\vadjust{\goodbreak} is referred to as the sparse case.
\item The parameter $\nu$
appearing here is often called degree of ill-posedness of the problem
(DIP). It appears here through conditions (\ref{blue}) and
(\ref{blue'}) which are essential in this problem. In
\cite{Kimoptimalspherical}, for instance, and very often in diverse
inverse problems, this DIP parameter is introduced with the help of the
eigenvalues of the operator (i.e., here the discrepancy of the
coefficients of $f_\eps$ in its expansion along the rotational
harmonics). In the following subsection, we prove that (\ref{blue})
and (\ref{blue'}) are in fact a consequence of (and even equivalent to)
the standard ``ordinary smooth'' condition.

\item The rates of convergence found here are standard in inverse
problems. They can be related to rates found in Kim and Koo \cite{Kimoptimalspherical} in
the same deconvolution problem, with a $\bL_2$ loss and constraints on
the spaces comparable to $B^{s}_{22}(M)$. In the deconvolution problem
on the interval, similar rates are found even for $\bL_p$ losses (with
standard modifications since the dimension here is 2 instead of 1):
see, for instance, Johnstone et al. \cite{jkpr}. These results
are proved to be minimax (see Kim and Koo
\cite{Kimoptimalspherical}) up to logarithmic factors, for the case
$p=2$ with a $B^{s}_{22}(M)$ constraint on the object to estimate. With
methods comparatively similar to those in Willer \cite{Thomas},
it could be proved that our rates are also minimax in the general case
(again up to logarithmic factors) if we assume condition (\ref{blue'})
(or, in other terms, that the DIP is exactly of the order $\nu$).

\item It is interesting to notice that Theorem \ref{upperinfty} proves
that the same algorithm is also working for the $\bL_\infty$ norm, with
a slightly more sophisticated proof. It is often the most useful loss
function in practice. The proof requires (\ref{blue'}) instead of
(\ref{blue}). We do not know if this condition is necessary.

\item It is worthwhile noticing that the procedure is adaptive, meaning
that it does not require a priori knowledge on the regularity (or
sparsity) of the function. It also adapts to nonhomogeneous smoothness
of the function. The logarithmic factor is a standard price to pay for
adaptation.

\item The procedure requires the knowledge of the constant $A$ which is
corresponding to the $\bL_\infty$ norm of the density $f_Z$. It is
obvious (because of the convolution) that $A\le M$. However, it should
be better to obtain a~procedure not depending on $M$ either. For that,
we advocate that $\|f_Z\|_\infty$ can be replaced by $\|\hat
f_Z(j_N))\|_\infty$ in practice where $\hat f_Z(j_N))$ is an
undersmoothed needlet estimator of the density $f_Z$ close to those
introduced in Baldi et al. \cite{BaldiKerkMariPica} but with
no thresholding and with the level $j_N$ chosen so that $2^{2j_N}
\simeq N/(\log N)^2$. Standard arguments similar as in Gin\'{e} and
Nickl~\cite{GineNickl10b}) imply that this random quantity
exponentially concentrates around $\|f_Z\|_\infty$. We can also adopt a
more straightforward strategy as detailed in Section \ref{sim-section}.
\end{enumerate}

\end{rem*}

\subsection{%
\texorpdfstring{Conditions (\protect\ref{blue}), (\protect\ref{blue'}) and the smoothness of {$f_\eps$}}%
{Conditions (37), (39) and the smoothness of {$f_\eps$}}}

Following Kim and Koo~(\cite{Kimoptimalspherical}, condition
(3.6)), we can define the smoothness of $f_{\eps}$ spectrally. We place
ourselves in the ``ordinary smooth'' case
%
\begin{equation}\label{smoothness}
\|(\hat{f}^{l}_{\eps})^{-1}\|_{\mathrm{op}}\leq d_0l^{\nu}\quad \mbox{and}\quad
\|\hat{f}^{l}_{\eps}\|_{\mathrm{op}}\leq d_1l^{-\nu}\qquad \mbox{as }
l\rightarrow\infty
\end{equation}
for some positive constants $d_0$, $d_1$ and nonnegative constant
$\nu$, and where the operator norm of the rotational Fourier transform
$\hat{f}^{l}_{\eps}$ is defined as
\[
\|\hat{f}^{l}_{\eps}\|_{\mathrm{op}}=\sup_{h\neq0, h\in\mathbb{H}_l}
\frac{\|\hat{f}^{l}_{\eps}h\|_2}{\|h\|_2},
\]
$\mathbb{H}_l$ being the $(2l+1)$-dimensional vector space spanned
by $\{Y^l_{m}\dvtx-l\leq m \leq l\}$.

The following proposition states that condition (\ref{blue}) [resp.
(\ref{blue'})] is verified
in the ordinary smooth case by the needlets system.

\begin{pro}\label{prop1}
If $\|(\hat{f}^{l}_{\eps})^{-1}\|_{\mathrm{op}}\le d_0l^{\nu}$, then there
exists a constant $C$ such that
\[
|\sigma_{j}|^{2}
:=A\sum_{ln}\biggl|\sum_m\psi_{j\eta}^{lm}\hat{f}^l_{\eps^{-1}mn}\biggr|^2\le C
2^{2j\nu}\qquad \forall j \ge0.
\]
If $\|\hat{f}^{l}_{\eps}\|_{op}\le d_0l^{-\nu}$, then there exists a
constant $C'$ such that
\[
|\sigma_{j}|^{2}
:=A\sum_{ln}\biggl|\sum_m\psi_{j\eta}^{lm}\hat{f}^l_{\eps^{-1}mn}\biggr|^2\ge C'
2^{2j\nu}\qquad \forall j \ge0.
\]
\end{pro}

The proof of this proposition is given in the supplement article
(Kerkyacharian et al. \cite{KerkPhamPicardA}).

Notice also that the super smooth case (corresponding to
exponential spectral decreasing) is also considered in Kim and Koo
\cite{Kimoptimalspherical}. We will not consider this case
here, although this could be done, basically because this case
corresponds to very poor rates of convergence (logarithmic in $N$). As
well, this case does not require a thresholding since the adaptation is
obtained almost for free.

We now give a brief review of some examples of smooth distributions
which are discussed in depth in Healy, Hendriks and Kim
\cite{HealyHendriksKim} and Kim and Koo
\cite{Kimoptimalspherical}.

\subsubsection{Rotational Laplace distribution}
This distribution can be viewed as an exact analogy on $\mathit{SO}(3)$ of the
Laplace distribution on $\mathbb{R}$. Spectrally, for some $\rho^2>0$,
this distribution is characterized by
%
\begin{equation}\label{Rotational_Laplace}
\hat{f}^l_{\eps,mn}=\bigl(1+\rho^2l(l+1)\bigr)^{-1}\delta_{mn}
\end{equation}
for $-l\leq m,n \leq l $ and $l=0,1,\ldots,$ and where $\delta_{mn}=1$ if
$m=n$ and $0$ otherwise.

\subsubsection{The Rosenthal distribution} %
This distribution has its origin in random walks in groups (for
details, see Rosenthal \cite{Rosenthal}).

If one considers the situation where $f_\eps$ is a $p$-fold convolution
product of conjugate invariant random for a fixed axis, then Rosenthal
(\cite{Rosenthal}, page 407) showed that
\[
\hat{f}^{l}_{\eps,mn}=\biggl(\frac{\sin(l+1/2)\theta}{(2l+1)\sin
\theta/2}\biggr)^{p}\delta_{mn}
\]
for $-l\leq m,n \leq l $ and $l=0,1,\ldots,$ and where $0<\theta\leq\pi$
and $p>0$.

\section{Practical performances}\label{sim-section}

In this section, we produce the results of numerical experiments on
the sphere $\mathbb{S}^{2}$. Numerical work has been conducted using
the spherical pixelization HEALPix software package. HEALPix provides
an approximate quadrature of the sphere with a number of data points of
order $C2^{2J}$ and a number of quadrature weights of order $\frac
{1}{C2^{2J}}$, for some positive constant $C$. This approximation is
considered as reliable enough and commonly used in astrophysics.

In the two examples below, we considered samples of cardinality
$N=1500$. The maximal resolution level $J$ is taken such that
$J=(1/2)\log_{2}(\frac{N}{\log N})$. In order not to have more
cubature points than observations, we set $J=3$ for $N=1500$. We recall
the expression of the estimate of the needlets coefficients of the
density of interest:
%
\begin{equation}
\hat{\beta}_{j\eta}=\frac{1}{N}\sqrt{\lambda_{j\eta}}\sum
_{l=2^{j-1}}^{2^{j+1}}b(l/2^j)\sum_{m=-l}^{l}\overline{Y^{l}_{m}(\xi
_{j\eta})}\sum_{n=-l}^{l}
\hat{f}^{l}_{\eps^{-1},mn}\sum_{u=1}^{N}{Y}^{l}_{n}(Z_{u}),
\end{equation}
where the quadrature weight are approximately uniform,
$\lambda_{j\eta}\simeq4\pi/(12.2^{2j})$. We replace the rotational
Fourier transform $(\hat{f}^{l}_{\eps})_{mn}:=\hat{f}^{l}_{\eps,mn}$
[defined in (\ref{rotationalFourierTransform})] by its empirical
version, more tractable in our simulation study. Note also that this
situation is very likely to occur for instance in the context of
astrophysics.

We precise again that $\hat{f}^{l}_{\eps^{-1},mn}$
denotes\vspace*{-2pt} the $(m,n)$ element of the matrix $(\hat{f}^{l}_{\eps})^{-1}:=
\hat{f}^{l}_{\eps^{-1}}$ which is the inverse\vspace*{-2pt} of the $(2l+1)\times
(2l+1)$ matrix $(\hat{f}^{l}_{\eps})$. In order to get the empirical\vspace*{-2pt}
version $\hat{f}^{l,N}_{\eps^{-1},mn}$ of $\hat{f}^{l}_{\eps^{-1},mn}$,
we have first to compute the empirical matrix $(\hat{f}^{l,N}_{\eps})$
then to inverse it to get the matrix
$(\hat{f}^{l,N}_{\eps})^{-1}:=\hat{f}^{l,N}_{\eps^{-1}}$. The $(m,n)$
entry of the matrix $(\hat{f}^{l,N}_{\eps})$ is given by the formula
\[
\hat{f}^{l,N}_{\eps,mn}=\frac{1}{N}\sum_{j=1}^{N}D^{l}_{m,n}(\eps_{j}),
\]
where the rotational harmonics $D^{l}_{m,n}$ have been defined in
(\ref{harmo_rotationnelle}). The $\eps_{j}$'s are i.i.d. realizations
of the variable $\eps\in \mathit{SO}(3)$.

For the generation of the random variable $\eps\in \mathit{SO}(3)$, we chose
Oz as the rotation axis and an angle $\phi$ following a uniform
law with different supports such as $[0,\pi/8]$, $[0,\pi/4]$,
$[0,\pi/2]$. The larger the support of distribution
is, the more intense the effect of the noise will be.

This particular choice of rotation matrix entails that in the
decomposition of an element of $\mathit{SO}(3)$ [see formula
(\ref{decompositionSO3})] the angles $\psi$ and $\theta$ are both equal
to zero. For this specific setting of perturbation, we deduce the
following form for the rotational harmonics:
\[
D^{l}_{m,n}(\eps_{j})=P^{l}_{mn}(1)e^{-ni\phi_{j}}=\delta
_{mn}e^{-in\phi_{j}},
\]
where $\phi_{j} \sim U[0,a]$ and $a$ is a positive constant which will
be specified later.

\textit{Choosing $\sigma_j$}.
As one may notice, the estimator $\hat{f}$ [see
(\ref{estimateur_de_f})] relies on the knowledge of $A$ which controls
the sup norm of the density $f_Z$ and appears in the formula of
$\sigma^2_j$, see (\ref{sigmaj}). Different ways to circumvent this
difficulty can be used, for instance, estimating it as explained above.
However, we have adopted in this section a different approach. The
quantity $\sigma_j^2$ constitutes\vspace*{-2pt} a~control of the variance of the
estimated coefficients $\hat{\beta}_{j\eta}$ as it is shown by
inequality (5) in the supplementary material (see Kerkyacharian
et al. \cite{KerkPhamPicardA}). Using this remark, instead of using an upper bound of
the variance of $\hat{\beta}_{j\eta}$, we will directly plug in an
estimation of this variance. Hence, $\sigma_j^2$ in (\ref{def-t-N}) is
replaced by
\[
v_{j\eta}^2=\frac{1}{(N-1)}\sum_{i=1}^{N}|G_{j\eta}(Z_i)-\widetilde
{G_{j\eta}}(Z)|^2,
\]
where
\[
G_{j\eta}(x)=\sum_{lm}\psi^{lm}_{j\eta}\sum_{n}\hat
{f}^{l}_{\varepsilon^{-1},mn}Y^l_{n}(x)
\]
and
\[
\widetilde{G_{j\eta}}(Z)=\frac{1}{N}\sum_{i=1}^{N}G_{j\eta}(Z_i).
\]
Remark that $\hat{\beta}_{j\eta}=\frac{1}{N}\sum_{i=1}^{N}
G_{j\eta}(Z_i)$.

Hence, for the reconstruction of the density $f$, we have the
following Needlet estimator:
\[
\hat{f}:=\frac{1}{|\mathbb{S}^{2}|}+\sum_{j=0}^{J}\sum_{\eta
=1}^{12.2^{2j}}\hat{\beta}_{j\eta}\psi_{j\eta}I_{\{|\hat{\beta
}_{j\eta}|\geq
\kappa t_{N}|v_{j\eta}|\}}.
\]

\textit{Choosing the tuning parameter $\kappa$.}
Before entering into the details of the numerical results, we will
dwell on the methodology used in this part to calibrate the tuning
paramater $\kappa$ in practice. In a first time, we focus on the
uniform density. Indeed, we make an upstream study with the uniform
density $f=\frac1 {4\pi} \mathbf{1}_{\mathbb{S}^2}$ which will
allow us to determine a reasonable value for $\kappa$ which will be
kept in the sequel for other types of densities. On the one hand, the
choice of the uniform density is not fortuitous because we know that in
theory the needlet coefficient $\langle f, \psi_{j\eta} \rangle=0$
for all $j$ and $\eta$. On the other hand, for an upstream study, a
prerequisite is to deal with a simple density, which is the case.
Hence, we test various values of $\kappa$ and see how many
coefficients survive thresholding. Then we look at the smallest value
of $\kappa$ for which all the estimated coefficients are killed. It
turns out that $\kappa=0.5$. Consequently for $\kappa=0.5,$ we have
noise-free reconstructions of the uniform density. Accordingly, this
particular value of $\kappa$ plays a benchmark value for the other
types of density that we highlight. In other words, in the Example 2,
which concerns the unimodal density we set $\kappa=0.5$.

We can form a parallel to what happens on the real line. Indeed, in a
certain way, this strategy for choosing $\kappa$ in practise meets up
with the universal threshold $\sqrt{2\log n}$ put forward by Donoho and
Johnstone  \cite{DonohoJohnstone} in the context of fixed
design regression on $\mathbb{R}$ in order to ``kill'' asymptotically all
the coefficients when estimating the zero function. Then, we compute
the $\mathbb{L}_2$ and $\mathbb{L}_\infty$ losses and give the graphic
reconstructions. This choice of $\kappa$ turns out to be very
reasonable and fruitful as the results prove to be good enough.

\begin{table}\tablewidth=265pt
\caption{Number of coefficients surviving thresholding\break for various
values of $\kappa$, $\phi\sim U[0,\pi/8]$} \label{Tableau}
\begin{tabular*}{265pt}{@{\extracolsep{\fill}}lcccc@{}}
\hline
& $\bolds{j=0}$ & $\bolds{j=1}$ & $\bolds{j=2}$ & $\bolds{j=3}$ \\
\hline
$\kappa=0.2$ & 0 & 7 & 30 & 110 \\
$\kappa=0.3$ & 0 &0 & 2 & 6 \\
$\kappa=0.4$ & 0 & 0 & 0 & 3 \\
\hline
\end{tabular*}
\end{table}
%

\begin{table}[b]
\tablewidth=265pt
\caption{Number of coefficients surviving thresholding for various
values of $\kappa$, $\phi\sim U[0,\pi]$}\label{tab2}
\begin{tabular*}{265pt}{@{\extracolsep{\fill}}lcccc@{}}
\hline
& $\bolds{j=0}$ & $\bolds{j=1}$ & $\bolds{j=2}$ & $\bolds{j=3}$ \\
\hline
$\kappa=0.2$ & 2 & 3 & 77 & 350 \\
$\kappa=0.3$ & 0 & 0 & 4 & 10 \\
$\kappa=0.4$ & 0 & 0 & 0 & 6 \\
\hline
\end{tabular*}
\end{table}

\begin{example}
In this first example, we consider the case of the
uniform density $f=\frac{1}{4\pi}\mathbf{1}_{\mathbb{S}^2}$. It is easy
to verify that $\beta_{j\eta}=\langle f,\psi_{j\eta}
\rangle_{\mathbb{L}_2}=0$ for every $j$ and every~$\eta$. Following
Baldi et al. \cite{BaldiKerkMariPica}, a simple way of assessing
the performance of the Needlet procedure is to count the number of
coefficients surviving thresholding.

We precise that both in the cases of an angle following a law
$U[0,\pi/8]$ or $U[0,\pi]$, all the coefficients are killed for
$\kappa=0.5$ (see Tables \ref{Tableau} and \ref{tab2}). Accordingly, we can conclude that the thresholding
procedure based on spherical needlets is very efficient.
\end{example}

\begin{example}
We will now deal with the example of a density of
the form $f(\omega)=ce^{-4|\omega-\omega_1|^2}$, with
$\omega_1=(0,1,0)$ and $c=1/0.7854$. The graph of $f$ in the spherical
coordinates $(\Phi,\Theta)$ ($\Phi={}$longitude, $0\leq\Phi\leq
2\pi$,
$\Theta={}$colatitude, $0\leq\Theta\leq\pi$) is given in Figure \ref{needlets}. We
also plot the noisy observations for different cases of perturbations.
For big rotation angles such that $\phi\sim U[0,\pi/4]$ or $\phi\sim
U[0,\pi/2]$, the observations tend to be spread over a large region on
the sphere and not to be concentrated in a specific region any more.
Consequently, denoising might prove to be difficult. In the context of
the deconvolution on the sphere, a~large amount of noise corresponds to
a rotation about the
Oz axis by a~large angle.

As motivated above, in the sequel, the tuning parameter $\kappa$ for
the Needlet procedure is set to $0.5$ both for computing the quadratic
loss, the $L_{\infty}$ loss and the graphic reconstructions.

%

\begin{table}[b]
\tablewidth=265pt
\caption{$\mathbb{L}_2$ and $\mathbb{L}_{\infty}$ estimated
losses}\label{tab3}
\begin{tabular*}{265pt}{@{\extracolsep{\fill}}lccc@{}}
\hline
& $\bolds{\phi\sim U[0,\pi/8]}$ & $\bolds{\phi\sim U[0,\pi/4]}$&
$\bolds{\phi\sim U[0,\pi/2]}$ \\
\hline
$\mathbb{L}_2$ & 0.3335 & 0.5523 & 0.7830 \\
$\mathbb{L}_{\infty}$& 0.1019 & 0.1677& 0.1928 \\
\hline
\end{tabular*}
\end{table}

\begin{figure}

\includegraphics{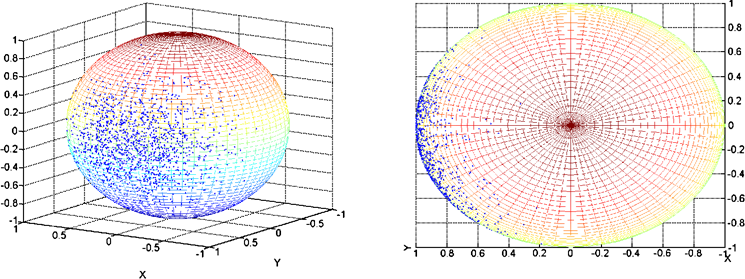}

\caption{Observations $\phi\sim U[0,\pi/8].$}\label{needlets}
\end{figure}

First of all, we have computed an estimate of the $\mathbb{L}_2$ and
the $\mathbb{L}_{\infty}$ norms of the difference between $\hat{f}$
our Needlet estimator and $f$ (see Table \ref{tab3}). For the quadratic loss, we took the
square root of the sum of the squares of the coefficients. As for the
$\mathbb{L}_{\infty}$ distance, we chose an almost uniform grid of
the sphere $\mathbb{S}^2$ given by the HEALPix pixelization program.
We recall that
\[
\|\hat{f}-f\|_{\infty}=\sup_{i=1,\ldots, L}|\hat{f}(\alpha
_i)-f(\alpha_i)|,
\]
where the $\alpha_i$ are the points of a uniform grid of the sphere.
Here, we chose $L=192=12.2^4$.
All the estimated losses were computed over $50$ runs and for $N=1500$
observations. We considered the three cases of noise level described above.
The results are summarized in the following table of errors which shows
that the estimator performs quite well. In particular, its performances
are deteriorating when the noise becomes more important (which was
expected) and give very good results in $\mathbb{L}_{\infty}$ norm.
We concentrated on particular phenomena instead of performing a large
scanning of the errors in very diverse situations because the
computations---although feasible---are rather costly in time when used
repeatedly. For instance, to our knowledge, and probably for the same
reason, no such study has been performed for the SVD method.
\end{example}

\textit{SVD versus needlet on a particular problem.}
A central issue in \textit{Astrophysics} is to detect the place of the peak of
the bell which in the present density case is localized in
$(\Theta=\pi/2, \Phi=\pi/2)$. For each case of noise,
we plot the observations both on the sphere and on the flattened sphere
and give the reconstructed density in the spherical coordinates for the
Needlet procedure and the SVD estimator (see Figures \ref{needlets}, \ref{fig4} and \ref{fig6}).

\begin{figure}

\includegraphics{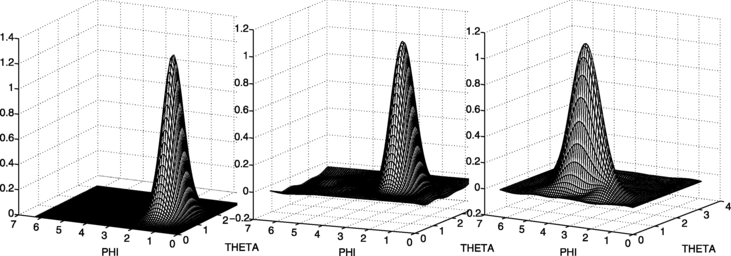}

\caption{The exact density, the Needlet procedure, the SVD method,
$\phi\sim U[0,\pi/8].$}\label{fig3}
\end{figure}

\begin{figure}[b]

\includegraphics{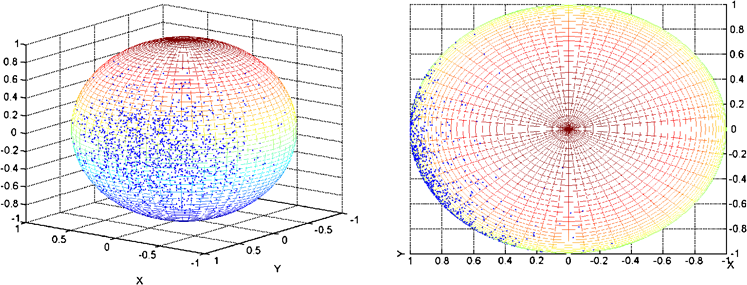}

\caption{Observations $\phi\sim U[0,\pi/4].$}\label{fig4}
\end{figure}

\begin{figure}

\includegraphics{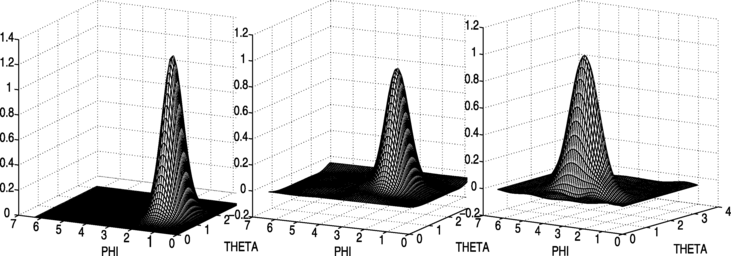}

\caption{The exact density, the Needlet procedure, the SVD method,
$\phi\sim U[0,\pi/4].$}\label{fig5}
\end{figure}

\begin{figure}[b]

\includegraphics{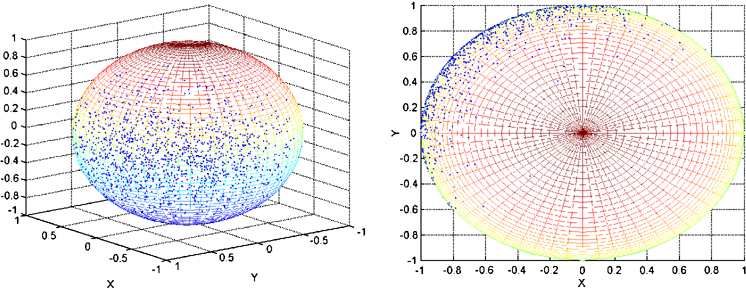}

\caption{Observations $\phi\sim U[0,\pi/2].$}\label{fig6}
\end{figure}

\begin{figure}

\includegraphics{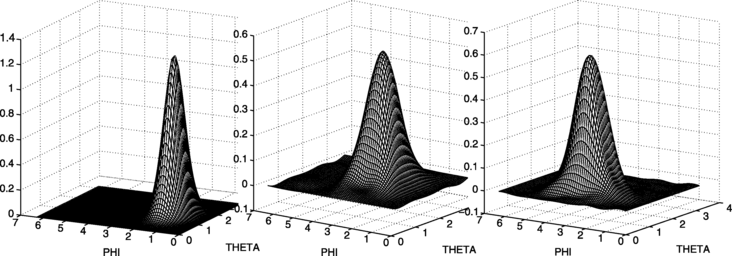}

\caption{The exact density, the Needlet procedure, the SVD method,
$\phi\sim U[0,\pi/2].$}\label{fig7}
\end{figure}

For each of the three groups of graphic reconstructions of the target
density presented below corresponding to three cases of noise, you will
find in order, the exact target density, then the one estimated with
the Needlet procedure and finally the density estimated with the SVD
method (see Figures \ref{fig3}, \ref{fig5} and \ref{fig7}).

%
%
At a closer inspection, we notice that the position of the peak of the
estimated bell is well localized by the Needlet procedure whatever the
amount of noise, it is especially true for $\phi\sim U[0,\pi/8]$. Only
in the case of the law $U[0,\pi/2]$, the longitude coordinate of the
peak tends to slightly move away from the true value. As for the SVD
method, for the three reconstructions, it fails to detect properly the
exact position of the peak. Therefore, even if in the case of big
rotations such that $\phi\sim U[0,\pi/4]$ and especially $\phi\sim
U[0,\pi/2]$, the Needlet procedure allows us to have a rather good
detection of the position of the peak. This is of course due to the
remarkable concentration of the needlet.
Of course, one remarks that the
base of the bell tends to become a bit larger when the noise increases,
this is due to the fact that the observations are not concentrated in a
specific region any longer, but the genuine form of the density is well
preserved.

\section{\texorpdfstring{Proof of Theorems \protect\ref{upperbound} and \protect\ref{upperinfty}}%
{Proof of Theorems 2 and 3}}
\label{proof}

In this proof, $C$ will denote an absolute constant which may change
from one line to the other.

%

We begin with the following proposition.
\begin{pro}\label{prop2}
For any $q\ge1$
there exist constants $ C,s_q, s'_q $, such that, as soon as
$2^{j}\le[\frac N{\log N}]^{1/2}$,
\begin{eqnarray}\label{bern}
\hspace*{45pt}\bP\{|\hat\beta_{j\eta}-\bjk|\ge\sigma_jv\}&\le&
2\exp\biggl\{-\frac{Nv^2}{2(1 +2/(\sqrt{3\pi} A) v2^{j})}\biggr\}\qquad
\forall v>0,\\\label{moment}
\bE|\hat\beta_{j\eta}-\bjk|^q &\le& s_q\biggl[\frac{\sigma
_{j}^2}N\biggr]^{q/2} ,\\\label{momsup}
\bE\sup_{\eta}|\hat\beta_{j\eta}-\bjk|^q &\le& s_q'(j+1)^q\biggl[\frac
{\sigma_{j}^2}N\biggr]^{q/2},\\\label{exp}
\bP(|\hat\beta_{j\eta}-\bjk|\ge\sigma_j{\kappa} t_N)&\le&
2N^{-\sqrt{3\pi}A\kappa/4}\qquad\forall \kappa\ge\sqrt{3\pi}A.
\end{eqnarray}
\end{pro}

The proof of this proposition is given in the supplementary material (see
Kerkyacharian et al. \cite{KerkPhamPicardA}).

\subsection{\texorpdfstring{Proof of Theorem \protect\ref{upperbound}}{Proof of Theorem 2}}

Now, to get the result of Theorem \ref{upperbound}, we begin by the
following decomposition:
\begin{eqnarray*}
\hspace*{2pt}\EE\|\hat f-f\|_p^p&\leq&
2^{p-1}\Biggl\{\bE\Biggl\|\sum_{j=-1}^{{J}}\sum_{\eta\in\mathscr Z_j}
\bigl(t(\hat\beta_{j\eta})-\bjk\bigr)\psijk\Biggr\|_p^p+\biggl\|\sum_{j>{J}}\sum_{\eta
\in
\mathscr Z_j} \bjk\psijk\biggr\|_p^p\Biggr\}\\
&=:& I+\mathit{II}.
\end{eqnarray*}
The term $\mathit{II}$ is easy to analyze, as follows. We observe first that
since $B^s_{\pi,r}(M) \subset B^s_{p,r}(M')$ for
$\pi\geq p$, this case will be assimilated to the case $\pi=p$ and
from now on, we will only consider $\pi\leq p$. Since $f$ belongs to
$B^s_{\pi,r}(M)$, using the embedding results recalled above in
(\ref{embedding}), we have that $f$ also belongs to $B^{s-(
2/\pi-2/p)}_{p,r}(M')$, for some constant $M'$ and for $\pi\leq
p$. Hence,\vspace*{-2pt}
\[
\biggl\|\sum_{j>{J}}\sum_{\eta\in\mathscr Z_j} \bjk\psijk\biggr\|_p\le C
2^{-J[s-2(1/\pi-1/p)]}.
\]
Then we only need to verify that $\frac{s-2(
1/\pi-1/p)}{\nu+1}$ is always larger that $ \mu$, which is not
difficult.

Indeed, on the first\vspace*{-1pt} zone $s\geq(\nu+1)(p/\pi-1)$. So,
$s+\nu+1 \geq(\nu+1)\frac{p}{\pi}$ which entails that
$\frac{s}{(s+\nu+1)}\leq\frac{s}{(\nu+1)p/\pi}$. We need to
check that $s-2(\frac1\pi-\frac1p) \geq\frac{s\pi}{p}$. We have that
$s-\frac2\pi+\frac2p-\frac{s\pi}{p}=2(\frac{s\pi}{2}-1)(\frac
1\pi
-\frac1p)\geq0$ since\vspace*{-2pt} $s\geq\frac2\pi$ and $p\geq\pi$.\\
On the second zone, we obviously have that $\frac{s-2(
1/\pi-1/p)}{\nu+1}$ is always larger than $ \mu=\frac{s-2/\pi
+2/p}{s+\nu-2/\pi+1}$.

Bounding the term I is more involved. Using the
triangular inequality together with H{\"{o}}lder inequality, and
property (\ref{item3}) for the second line, we get\vspace*{-2pt}
\begin{eqnarray*}
I&\le&2^{p-1}J^{p-1} \sum_{j=-1}^{{J}}\bE\biggl\|\sum_{\eta\in\mathscr
Z_j} \bigl(t(\hat\beta_{j\eta})-\bjk\bigr)\psijk\biggr\|_p^p\\
&\le&2^{p-1}J^{p-1}C\sum_{j=-1}^{{J}}\sum_{\eta\in\mathscr Z_j}\bE
|t(\hat\beta_{j\eta})-\bjk|^p\|\psijk\|_p^p.
\end{eqnarray*}
Now, we separate four cases:
\begin{eqnarray*}
&&\sum_{j=-1}^{{J}}\sum_{\eta\in\mathscr Z_j}\bE
|t(\hat\beta_{j\eta})-\bjk|^p\|\psijk\|_p^p\\
&&\qquad=\sum_{j=-1}^{{J}}\sum_{\eta\in\mathscr Z_j} \bE
|t(\hat\beta_{j\eta})-\bjk|^p\|\psijk\|_p^p
\bigl\{I\{|\hat\beta_{j\eta}|\ge{\kappa} t_N{\sigma_j}\}\\
&&\hphantom{\qquad=\sum_{j=-1}^{{J}}\sum_{\eta\in\mathscr Z_j} \bE
|t(\hat\beta_{j\eta})-\bjk|^p\|\psijk\|_p^p
\bigl\{}{}+I\{|\hat\beta_{j\eta}|< {\kappa} t_N{\sigma_j}\}\bigr\}\\
&&\qquad\le \sum_{j=-1}^{{J}}\sum_{\eta\in\mathscr Z_j}\biggl[\bE
|\hat\beta_{j\eta}-\bjk|^p\|\psijk\|_p^p I\{|\hat\beta_{j\eta}|\ge
{\kappa} t_N{\sigma_j}\}\\
&&\hphantom{\qquad\le \sum_{j=-1}^{{J}}\sum_{\eta\in\mathscr Z_j}\biggl[}{}\times \biggl\{ I\biggl\{|\bjk|\ge\frac\kappa2 t_N{\sigma_j}\biggr\}+ I\biggl\{|\bjk|< \frac\kappa2 t_N{\sigma_j}\biggr\}\biggr\}\\
&&\hphantom{\qquad\le \sum_{j=-1}^{{J}}\sum_{\eta\in\mathscr Z_j}\biggl[}{}+ |\bjk|^p\|\psijk\|_p^pI\{|\hat\beta_{j\eta}|< {\kappa}
t_N{\sigma_j}\}\bigl\{I\{|\bjk|\ge2{\kappa} t_N{\sigma_j}\}\\
&&\hphantom{\hphantom{\qquad\le \sum_{j=-1}^{{J}}\sum_{\eta\in\mathscr Z_j}\biggl[}{}+ |\bjk|^p\|\psijk\|_p^pI\{|\hat\beta_{j\eta}|< {\kappa}
t_N{\sigma_j}\}\bigl\{}{}+I\{|\bjk|< 2{\kappa} t_N{\sigma_j}\}\bigr\}\biggr]\\
&&\qquad= :\mathit{Bb}+\mathit{Bs}+\mathit{Sb}+\mathit{Ss} .
\end{eqnarray*}

We have the following upper bounds for the terms $\mathit{Bs}$ and $\mathit{Sb}$:
\begin{eqnarray*}
\mathit{Bs}&\leq& CN^{-\sqrt{3\pi}A\kappa/16},\\
\mathit{Sb} &\leq& CN^{-\sqrt{3\pi}A\kappa/2+p/(\nu+1)}.
\end{eqnarray*}

It is easy to check that in any cases if $
\sqrt{3\pi}A\kappa>\max\{8p, 2p+1\}$ the terms $\mathit{Bs}$ and $\mathit{Sb}$ are
smaller than the rates announced in the theorem. For the details of the
above upper bounds of the terms $\mathit{Bs}$ and $\mathit{Sb}$, see the
supplementary
material
(Kerkyacharian et al. \cite{KerkPhamPicardA}).

Using Proposition \ref{prop2}, we have
\begin{eqnarray*}
\mathit{Bb}&\le&
C\sum_{j=-1}^{{J}}\sum_{\eta\in\mathscr Z_j}\sigma_j^p
N^{-p/2}\|\psijk\|_p^p I\biggl\{|\bjk|\ge\frac\kappa2 t_N{\sigma_j}\biggr\},\\
\mathit{Ss}&\le&\sum_{j=-1}^{{J}}\sum_{\eta\in\mathscr
Z_j}|\bjk|^p\|\psijk\|_p^pI\{|\bjk|< 2{\kappa} t_N{\sigma_j}\}.
\end{eqnarray*}

Using again Proposition \ref{prop2}, (\ref{lp}) and condition
(\ref{blue}) for any $p\ge z\ge0$:
\begin{eqnarray*}
\mathit{Bb}&\leq& C N^{-p/2}\sum_{j=-1}^{{J}}\sigma_j^p2^{j(p-2)}\sum_{\eta
\in
\mathscr Z_j} I\biggl\{|\bjk|\ge\frac\kappa2 t_N{\sigma_j}\biggr\}\\
&\leq& CN^{-p/2}\sum_{j=-1}^{{J}}\sigma_j^p2^{j(p-2)}\sum_{\eta\in
\mathscr Z_j}|\bjk|^z [ t_N{\sigma_j}]^{-z} \\
&\leq& C\tep^{p-z}\sum_{j=-1}^{{J}}2^{j[\nu(p -z)+p-2]}\sum_{\eta\in
\mathscr
Z_j}|\bjk|^z.
\end{eqnarray*}
Also, for any $p\geq z\geq0$,
\begin{eqnarray*}
\mathit{Ss} &\leq& C\sum_{j=-1}^{{J}}2^{j(p-2)}\sum_{\eta\in\mathscr
Z_j}|\bjk|^z \sigma_j^{p-z}[ t_N]^{p-z}\\
&\leq&
C [ t_N]^{p-z}\sum_{j=-1}^{{J}}2^{j(\nu(p -z)+p-2)}\sum_{\eta\in
\mathscr Z_j}|\bjk|^z.
\end{eqnarray*}

So in both cases we have the same bound to investigate. We will write
this bound on the following form (forgetting the constant):
\begin{eqnarray*}
A+B&:=&\tep^{p-z_1}\Biggl[\sum_{j=-1}^{{j_0}}2^{j[\nu(p -z_1)+p-2]}\sum
_{\eta\in\mathscr Z_j}|\bjk|^{z_1}\Biggr]\\
&&{}+\tep^{p-z_2}\Biggl[\sum
_{j=j_0+1}^{{J}}2^{j[\nu(p -z_2)+p-2]}\sum_{\eta\in\mathscr
Z_j}|\bjk|^{z_2}\Biggr].
\end{eqnarray*}

The constants $z_i$ and $j_0$ will be chosen depending on the cases,
with the only constraint $p\geq z_i\geq0$.

We recall that we only need to investigate the case $p\ge\pi$, since
when $p\le\pi$, $B^s_{\pi r}(M)\subset B^s_{p r}(M')$.

Let us first consider the case where $s\ge(\nu+1)(\frac{p}{\pi}-1)$,
put
\[
q=\frac{p(\nu+1)}{s+\nu+1},
\]
and observe that on the considered
domain, $q\le\pi$ and $p>q$. In the sequel, it will be useful to
observe that we have $s=(\nu+1)(\frac{p}q-1)$. Now, taking
$z_2=\pi$, we get
\[
B\le\tep^{p-\pi}\Biggl[\sum_{j=j_0+1}^{{J}}2^{j[\nu(p -\pi)+p-2]}\sum
_{\eta\in\mathscr Z_j}|\bjk|^{\pi}\Biggr].
\]
Now, as
\[
\frac{p}{q}-\frac2\pi+\nu\biggl(\frac pq-1\biggr)=s+1-\frac2\pi
\]
and
\[
\sum_{\eta\in\mathscr Z_j}|\bjk|^{\pi}=2^{-j\pi(s+1-2/\pi
)}\tau_j^{\pi}\label{bspi},
\]
with $(\tau_j)_j \in l_r$ (this last thing is a consequence of the
fact that $f\in B^s_{\pi,r}(M)$), we can write
\begin{eqnarray*}
B &\leq&\tep^{p-\pi}\sum_{j=j_0+1}^{J}2^{j p(1-\pi/
q)(\nu+1)}\tau_j^\pi\\
&\leq& C\tep^{p-\pi} 2^{j_0 p(1-\pi/q)(\nu+1)}.
\end{eqnarray*}
The last inequality is true for any $r\ge1$ if $\pi> q$ and for $r\le
\pi$ if $\pi=q$. Notice that $\pi=q$ is equivalent to $ s = (\nu
+1)(\frac{p}{\pi} - 1). $ Now if we choose $j_0$ such that $2^{j_0
(p/q) (\nu+1)}\sim\tep^{-1}$, we get the bound
\[
\tep^{p-q},
\]
which exactly gives the rate announced in the theorem for this case. As
for the first part of the sum (before $j_0$), we have, taking now
$z_1=\tilde q$, with $\tilde q\le\pi$ (and also $\tilde q\le p$ since
we investigate the case $p\ge\pi$), so that\break
$[\frac1{2^{2j}}\sum_{\eta\in\mathscr Z_j}|\bjk|^{\tilde q}]^{1/\tilde q}
\le[\frac1{2^{2j}}\sum_{\eta\in\mathscr Z_j}|\bjk|^\pi]^{1/\pi}$, we get
\begin{eqnarray*}
A&\leq& \tep^{p-\tilde q}\Biggl[\sum_{j=-1}^{{j_0}}2^{j[\nu(p -\tilde
q)+p-2]}\sum_{\eta\in\mathscr Z_j}|\bjk|^{\tilde q}\Biggr]\\
&\leq&\tep^{p-\tilde q}\Biggl[\sum_{j=-1}^{{j_0}}2^{j[\nu(p -\tilde
q)+p-2\tilde q/\pi]}\biggl[\sum_{\eta\in\mathscr
Z_j}|\bjk|^{\pi}\biggr]^{\tilde q/\pi}\Biggr]\\
&\leq&\tep^{p-\tilde q}\sum_{j=-1}^{{j_0}}2^{j[(\nu+1) p(1-
\tilde q/q)]} \tau_j^{\tilde q}\\
&\leq& C\tep^{p-\tilde q}2^{j_0p[(\nu+1) (1-\tilde q/q)]}\\
&\leq& C\tep^{p-q}.
\end{eqnarray*}
The last two lines are valid if $\tilde q$ is chosen strictly smaller
than $q$ (this is possible since $\pi\ge q$).

Let us now consider the case where $s< (\nu+1)(\frac{p}{\pi}
- 1)$, and choose now
\[
q=p\frac{\nu+1 -2/p}{s+\nu-2/\pi+1}.
\]
In such a way that we\vspace*{-1pt} easily
verify that $p-q=p\frac{s-2/\pi+2/p}{1+\nu+s-2/\pi},q-\pi
=\break \frac{(p-\pi)(1+\nu)-\pi s}{s+\nu-2/\pi+1}> 0$. Furthermore,\vspace*{1pt} we
also have $ s+1-\frac2\pi=\frac{ p}{q}-\frac2q+\nu(\frac
{p}q-1)$.
Hence, taking $z_1=\pi$ and using again the fact that $f$ belongs to $
B^s_{\pi,r}(M)$,
\begin{eqnarray*}
A&\le&\tep^{p-\pi}\Biggl[\sum_{-1}^{{j_0}}2^{j[\nu(p -\pi)+p-2]}\sum
_{\eta\in\mathscr Z_j}|\bjk|^{\pi}\Biggr]\\
&\le&\tep^{p-\pi}\sum_{-1}^{{j_0}}2^{j[(\nu+1-2/p)
(p/q)(q-\pi)]}\tau_j^\pi\\
&\le& C \tep^{p-\pi}2^{j_0[(\nu+1-2/p)(p/q)(q-\pi)]}.
\end{eqnarray*}
This is true since $\nu+1-\frac2p$ is also strictly positive since
$\nu+1>\frac s{p/\pi-1}\ge\frac2{ p-\pi}\ge\frac2p$. If we now
take $2^{j_0 (p/q) (\nu+1-2/p)}\sim\tep^{-1}$, we get
the bound
\[
\tep^{p-q},
\]
which is the rate announced in the theorem for this case.\vspace*{1pt}

Again, for $B$, we have, taking now $z_2=\tilde q>q(>\pi)$\vspace*{2pt}
\[
B\le\tep^{p-\tilde q}\Biggl[\sum_{j=j_0+1}^{{J}}2^{j[\nu(p -\tilde
q)+p-2]}\sum_{\eta\in\mathscr Z_j}|\bjk|^{\tilde
q}\Biggr].
\]
But
\[
\sum_{\eta\in\mathscr Z_j}|\bjk|^{\tilde q}\le C2^{-j\tilde{q}(s+1
- 2/\pi)}\tau^{\tilde{q}}_j,
\]
and $ s+1-\frac2\pi=\frac{ p}{q}-\frac
2q+\nu(\frac{p}q-1)$, hence
\begin{eqnarray*}
B &\le& C\tep^{p-\tilde q}\sum_{j=j_0+1}2^{j[(\nu+1-2/p)
(p/q)(q-\tilde q)]}\tau_j^{{\tilde q}}
\\
&\le& C\tep^{p-\tilde q}2^{j_0[(\nu+1-2/p)(p/q)(q-\tilde q)]}
\\
&\le& C\tep^{p-q},
\end{eqnarray*}
which completes the proof of Theorem \ref{upperbound}.
\subsection{\texorpdfstring{Proof of Theorem \protect\ref{upperinfty}}{Proof of Theorem 3}}

The proof of this theorem is entirely given in the supplementary material (see
Kerkyacharian et al. \cite{KerkPhamPicardA}).

\begin{supplement}[id=suppA]
\stitle{Supplement to ``Localized spherical deconvolution''\\}
\slink[doi]{10.1214/10-AOS858SUPP}
\sdatatype{.pdf}
\sfilename{supplementA\_Annals.pdf}
\sdescription{We give in the supplement some technical details
for the understanding of the proofs of Propositions \ref{prop1} and \ref{prop2} and of
Theorems \ref{upperbound} and \ref{upperinfty}.}
\end{supplement}

\section*{Acknowledgments}

The authors would like to thank Erwan Le Pennec for many helpful
discussions and suggestions concerning the simulations. We would also
like to thank an Associate Editor and two anonymous referees for their
insightful comments on a first draft of this work.

%


\printaddresses

\end{document}